\documentstyle[12pt]{article}

\newcommand{\sect}[1]{\setcounter{equation}{0}\section{#1}}

\def\be{\begin{equation}}
\def\ee{\end{equation}}
\def\bea{\begin{eqnarray}}
\def\eea{\end{eqnarray}}

\def\R{{\rm I\kern-.2em R}}

\def\>#1{{\bf #1}}                 

\def\1{\'{\i}}                           

\def\back{\!\!\!\!\!\!}

\begin{document}

\thispagestyle{empty}

 \hfill \ 
\bigskip
\bigskip
\bigskip
\bigskip
\bigskip

\begin{center} 

{\LARGE{\bf{Non-standard  quantum $so(3,2)$}}}

{\LARGE{\bf{and its contractions}}} 

\end{center}

\bigskip\bigskip

\begin{center}  
Francisco J. Herranz 
\end{center}

\begin{center} {\it Departamento de F\1sica\\ E.U. Polit\'ecnica,  
Universidad
de Burgos\\ 
E-09006, Burgos, Spain}

{\it  e-mail: fteorica@cpd.uva.es}

\end{center}

\bigskip\bigskip\bigskip

\begin{abstract} 
A full (triangular) quantum deformation of $so(3,2)$ is presented by
considering this algebra as the conformal algebra of the 2+1
dimensional Minkows\-kian spacetime. Non-relativistic  contractions are
analysed and used to obtain quantum Hopf structures for the conformal
algebras of the 2+1 Galilean and Carroll spacetimes. Relations between
the latter   and the null-plane quantum Poincar\'e algebra are studied.
\end{abstract}

\newpage


\sect {Introduction}

The Lie algebra ${\cal M}_3$ of the group of conformal transformations in 
the 2+1 Minkows\-kian spacetime is a ten-dimensional Lie algebra isomorphic
to $so(3,2)$. We consider the basis
$\{J,P_0,P_i,K_i,C_0,C_i,D\}$ $(i=1,2)$ where $J$ generates rotations,
$P_0$ time translations, $P_i$ space translations, $K_i$ boosts, $C_0$
and $C_i$ special conformal transformations, and $D$ dilations.
The  Lie brackets of ${\cal M}_3$ are
 \be
\begin{array}{lll}
[J,K_i]=\epsilon_{ij}K_j&\qquad [J,P_i]=\epsilon_{ij}P_j&\qquad
[J,C_i]=\epsilon_{ij}C_j\cr
[K_i,P_0]=P_i&\qquad [K_i,P_j]=\delta_{ij} P_0&\qquad
[K_1,K_2]=-J\cr
[K_i,C_0]=C_i&\qquad [K_i,C_j]=\delta_{ij} C_0&\qquad
[P_0,C_0]=D\cr
[P_0,C_i]=-K_i&\qquad [C_0,P_i]=-K_i&\qquad
[P_i,C_j]=-\delta_{ij}D+\epsilon_{ij} J\cr
[D,P_\mu]=P_\mu&\qquad [D,C_\mu]=-C_\mu&\qquad
[P_\mu,P_\nu]=0\qquad [C_\mu,C_\nu]=0\cr
[J,P_0]=0&\qquad [J,C_0]=0&\qquad [D,J]=0\qquad\quad\! [D,K_i]=0
\end{array}
\label{aa}
\ee
where $\epsilon_{ij}$ is antisymmetric with $\epsilon_{12}=1$, 
$\epsilon_{21}=-1$, $\epsilon_{ii}=0$, and from now on we assume
that $\mu,\nu=0,1,2$ and  $i,j=1,2$. The 2+1 Poincar\'e algebra,
${\cal P}(2+1)\equiv \langle J,P_0,P_i,K_i\rangle$, is a Lie subalgebra of
${\cal M}_3$; moreover, if we add the dilation generator $D$ to ${\cal
P}(2+1)$ we get the Weyl Lie subalgebra  $\overline{{\cal P}}(2+1)$.
Hence we have the sequence ${\cal P}(2+1)\subset \overline{{\cal P}}(2+1)
\subset {\cal M}_3$.

The two non-relativistic limits of the Poincar\'e algebra ${\cal P}(2+1)$ are
the Galilean ${\cal G}(2+1)$ and Carroll ${\cal C}(2+1)$ algebras which 
correspond, in this order,  to a speed-space and a speed-time contraction of
${\cal P}(2+1)$  \cite{BLL}. These contraction processes can be
implemented at a conformal level in order to obtain the conformal
algebras of the 2+1 Galiean and  Carroll spacetimes \cite{conform}, here
denoted
${\cal G}_3$ and ${\cal C}_3$, by applying the following mappings to the
generators of ${\cal M}_3$:
 \be
\begin{array}{lllll}
\mbox{Speed-space contr.:}&\   J\to J&\  P_0\to P_0&\ 
C_0\to C_0&\  D\to D  \cr
{\cal M}_3\to {\cal G}_3 &\  P_i\to \varepsilon P_i&\  K_i\to 
 \varepsilon K_i&
\  C_i\to \varepsilon C_i &  
\end{array}
\label{ab}
\ee
 \be
\begin{array}{lllll}
\mbox{Speed-time contr.:}&\   J\to J&\  P_i\to  P_i&\ 
C_i\to  C_i&\  D\to D  \cr
{\cal M}_3\to {\cal C}_3 &\  P_0\to \varepsilon P_0&\  K_i\to 
 \varepsilon K_i&
\  C_0\to \varepsilon C_0 &  
\end{array}
\label{ac}
\ee
Once these transformations have been performed on  the Lie brackets
(\ref{aa}) we get after the limit $\varepsilon\to 0$ the commutation
relations of 
  ${\cal G}_3$ and ${\cal C}_3$. For the Galilean case, the non-vanishing
commutators are:
 \be
\begin{array}{lll}
[J,K_i]=\epsilon_{ij}K_j&\quad [J,P_i]=\epsilon_{ij}P_j&\quad
[J,C_i]=\epsilon_{ij}C_j\cr
[K_i,P_0]=P_i&\quad  [K_i,C_0]=C_i&\quad  
[P_0,C_0]=D\cr
[P_0,C_i]=-K_i&\quad [C_0,P_i]=-K_i&\quad
[D,P_\mu]=P_\mu \quad\ [D,C_\mu]=-C_\mu .
\end{array}
\label{ad}
\ee
The conformal Galilean Lie algebra ${\cal G}_3$ is isomorphic to  
$t_6(so(2)\oplus so(2,1))$ (the structure of this type of algebras is
described in  \cite{WB,coho}). We also have a sequence of subalgebras ${\cal
G}(2+1)\subset
\overline{{\cal G}}(2+1)
\subset {\cal G}_3$, where $\overline{{\cal G}}(2+1)$ is the 2+1 Galilean
algebra with dilation. 

Likewise, we obtain the conformal Carroll algebra
${\cal C}_3$  with non-zero Lie brackets given by:
 \be
\begin{array}{lll}
[J,K_i]=\epsilon_{ij}K_j&\quad [J,P_i]=\epsilon_{ij}P_j&\quad
[J,C_i]=\epsilon_{ij}C_j\cr
[K_i,P_i]=P_0&\quad [K_i,C_i]=C_0&\quad  
[P_0,C_i]=-K_i\quad\  [C_0,P_i]=-K_i\cr
[D,P_\mu]=P_\mu&\quad [D,C_\mu]=-C_\mu&\quad
[P_i,C_j]=-\delta_{ij}D+\epsilon_{ij} J  .
\end{array}
\label{ae}
\ee
The embedding  ${\cal
C}(2+1)\subset \overline{{\cal C}}(2+1) \subset {\cal C}_3$  is
easily verified  ($\overline{{\cal C}}(2+1)$ means the Carroll Weyl
subalgebra). The conformal  algebra ${\cal C}_3$ turns out to be isomorphic
to the 3+1 Poincar\'e algebra $iso(3,1)$. Recall that, in general, 
kinematical symmetries in $N+1$ dimensions can been seen as conformal
symmetries in $N$ dimensions
\cite{conform}.

Non-standard quantum deformations for these conformal algebras haven been
already obtained for the 1+1 case \cite{Beyond,av}  being inspired in 
the well known non-standard or Jordanian quantum $sl(2,\R)$ algebra
\cite{Demidov,Zakr,Ohn}. Their  underlying Lie bialgebras come from  
classical $r$-matrices which satisfy the classical Yang--Baxter equation.
The results so obtained show that non-standard deformations are naturally
adapted to a conformal basis, although for  the particular case of the
quantum Poincar\'e algebra an alternative interpretation  has been
considered in a null-plane framework \cite{Null}.

An analysis of non-standard conformal Lie bialgebras for higher
dimensions  can be found in \cite{Lukierski}  where the deformation
parameters are interpreted as fundamental mass parameters. However, to our
knowledge, no explicit non-standard quantum Hopf structure for the conformal
algebra further the 1+1 case ($so(2,2)$) has been given yet. In this letter
we solve this problem for a precise non-standard   quantum deformation of
${\cal M}_3$. To begin with we consider in section 2 the 2+1 conformal Lie
bialgebra which generalizes that introduced in  
\cite{Beyond,av}, and we study its possible non-relativistic Lie bialgebra
contractions.  It is shown that there is a unique possible (coboundary)
contraction for each conformal bialgebra ${\cal
G}_3$ and ${\cal C}_3$. In section 3 the quantum Hopf structure of ${\cal
M}_3$ is introduced and those corresponding to  ${\cal G}_3$ and ${\cal C}_3$
are obtained via contraction in section 4.  All  of them have as Hopf
subalgebra the corresponding kinematical algebra together the dilation
generator, but not the kinematical algebra itself (a feature already pointed
out in \cite{Dobrev}); hence only the Weyl subalgebra
is promoted to a Hopf subalgebra. The quantum conformal Carroll
algebra is related with the 3+1 null-plane quantum Poincar\'e algebra; this
fact allows to  get  its universal $R$-matrix from the results given in
\cite{rnull}.


\sect {Conformal Lie bialgebras}

 The classical $r$-matrices underlying the non-standard
quantum deformations of \break $sl(2,\R)\equiv\langle P_0,C_0,D\rangle$  and 
$so(2,2)\equiv\langle P_0,P_1,K_1,C_0,C_1,D\rangle$
\cite{Beyond,av} can be written  as
\be
r=z D\wedge P_0\qquad r=z(D\wedge P_0 + K_1\wedge P_1)
\label{ba}
\ee
where $z$ is the deformation parameter.  The generalization of these
expressions for ${\cal M}_3\simeq so(3,2)$ reads
\be
r=z(D\wedge P_0+K_1\wedge P_1 + K_2\wedge P_2 + J\wedge P_2)
\label{bb}
\ee
which fulfills the classical Yang--Baxter equation (the presence of the term
$J\wedge P_2$ is essential for this purpose). The cocommutator of  a
generator $X$ is obtained as $\delta(X)=[1\otimes X + X\otimes 1,r]$, namely,
\be
\begin{array}{l}
\delta(P_0)=0\qquad \delta(P_1)=0\cr
\delta(P_2)=-zP_2\wedge P_1\qquad \delta(J)=-zJ\wedge P_1\cr
\delta(D)=z(D\wedge P_0+K_1\wedge P_1 + K_2\wedge P_2 + J\wedge P_2)\cr
\delta(K_1)=z(K_1\wedge P_0+D\wedge P_1 - K_2\wedge P_2 - J\wedge P_2)\cr
\delta(K_2)=z(K_2\wedge P_0+J\wedge P_0 + J\wedge P_1 + K_1\wedge P_2
+D\wedge P_2)\cr
\delta(C_0)=z(C_0\wedge P_0-C_1\wedge P_1 -C_2\wedge P_2 - J\wedge K_2)\cr
\delta(C_1)=z(C_1\wedge P_0-C_0\wedge P_1 -C_2\wedge P_2 - J\wedge K_2)\cr
\delta(C_2)=z(C_2\wedge P_0+C_1\wedge P_2 -C_0\wedge P_2 + J\wedge K_1
+J\wedge D) .
\end{array}
\label{bc}
\ee

In order to analyse the possible non-relativistic contractions of this
conformal Lie bialgebra one has to consider the  Lie algebra transformations
(\ref{ab}) and (\ref{ac}) together a mapping on the deformation parameter:
$z\to \varepsilon^{-n} z$ where $n$ is any real number \cite{LBC}.  The
result is that there exists a unique minimal value $n_0$ of $n$ for each
contraction from ${\cal M}_3$ to  ${\cal G}_3$ and ${\cal C}_3$  in such way
the classical $r$-matrix and the cocommutators do not present divergencies:
\bea
&&{\cal M}_3\to {\cal G}_3:\quad z\to \varepsilon^{-2} z\qquad
(n_0=2)\label{bbd}\\ 
&&{\cal M}_3\to {\cal C}_3:\quad z\to \varepsilon^{-1} z \qquad
(n_0=1).
\label{bd}
\eea
For   $n>n_0$  the contracted $r$-matrix and  
cocommutators go to zero and for $n<n_0$ they diverge.

The classical (non-standard) $r$-matrix and cocommutators of the conformal
Gali\-lean Lie bialgebra so obtained are given by:
\be
r=z(K_1\wedge P_1 + K_2\wedge P_2 )
\label{be}
\ee
\be
\begin{array}{l}
\delta(X)=0\qquad \mbox{for}\ X\in\{J,P_\mu,K_i,C_i\}\cr
\delta(D)=z(K_1\wedge P_1 + K_2\wedge P_2 )\cr
\delta(C_0)=-z(C_1\wedge P_1 + C_2\wedge P_2 ) ,
\end{array}
\label{bf}
\ee
and for the Carroll case we get:
\be
r=z(D\wedge P_0+K_1\wedge P_1 + K_2\wedge P_2)
\label{bg}
\ee
\be
\begin{array}{l}
\delta(X)=0\qquad \mbox{for}\ X\in\{J,P_\mu\}\cr
\delta(Y)=z \, Y\wedge P_0\qquad \mbox{for}\ Y\in\{K_i,C_0\}\cr
\delta(D)=z(D\wedge P_0+K_1\wedge P_1 + K_2\wedge P_2)\cr
\delta(C_1)=z(C_1\wedge P_0-C_0\wedge P_1  - J\wedge K_2)\cr
\delta(C_2)=z(C_2\wedge P_0 -C_0\wedge P_2 + J\wedge K_1) .
\end{array}
\label{bh}
\ee

It is worth remarking that in each of the above Lie bialgebras the
corresponding Weyl subalgebra $\{J,P_\mu,K_i,D\}$ is preserved  at a
bialgebra level. Note also that the cocommutator of $D$ coincides with the
classical
$r$-matrix.


\sect {Quantum conformal Hopf algebra}

We proceed to introduce the Jordanian quantum deformation of the conformal
Minkowskian bialgebra, $U_z({\cal M}_3)$, in two steps. Firstly we close the
Hopf structure for the Weyl subalgebra, and secondly we complete the quantum
deformation with the expressions involving the special conformal
transformations. None direct procedure as the deformation embedding method
(applied for instance to the null-plane  quantum Poincar\'e algebra
\cite{Null}) seems to be useful now, so that we are forced to deduce
formerly the coproduct
$\Delta$ by solving the coassociativity condition
\be
(1\otimes \Delta)\Delta=(\Delta\otimes 1)\Delta
\ee
and by taking into account that
the cocommutators (\ref{bc}) are related with the first order of $\Delta$ on
$z$, $\Delta_{(1)}$, by means of 
\be
\delta=(\Delta_{(1)}-\sigma \circ \Delta_{(1)})\qquad
\mbox{where}\quad  \sigma(a\otimes b)=b\otimes a.
\ee
 Afterwards, the deformed commutation rules
follow by imposing the coproduct to be an algebra homomorphism, this is,
$\Delta([X,Y])=[\Delta(X),\Delta(Y)]$. 

In the sequel we write down the
coproduct and the commutation relations for 
$U_z({\cal M}_3)$; the counit is  trivial  and the antipode can be easily
derived from these results so we omit them.

\noindent
a) Weyl Hopf subalgebra $U_z(\overline{{\cal P}}(2+1))$.
\bea
&&\Delta(P_0)=1\otimes P_0 + P_0\otimes 1\qquad 
\Delta(P_1)=1\otimes P_1 + P_1\otimes 1\cr
&&\Delta(P_2)=1\otimes P_2 + P_2\otimes e^{-zP_1}\qquad 
\Delta(J)=1\otimes J + J\otimes e^{-zP_1}\cr
&&\Delta(D)=1\otimes D + D\otimes e^{zP_0}\cosh zP_1  + K_1\otimes
e^{zP_0}\sinh zP_1 \cr
&&\qquad\qquad +z(J+K_2)\otimes e^{zP_0}P_2 +\frac {z^2}2(K_1+D)\otimes 
e^{zP_0} e^{zP_1}P_2^2 \label{ca}\\
&&\Delta(K_1)=1\otimes K_1  + K_1\otimes e^{zP_0}\cosh zP_1  + D\otimes
e^{zP_0}\sinh zP_1 \cr
&&\qquad\qquad -z(J+K_2)\otimes e^{zP_0}P_2 -\frac {z^2}2(K_1+D)\otimes 
e^{zP_0} e^{zP_1}P_2^2\cr
&&\Delta(K_2)=1\otimes K_2  + (J+K_2)\otimes e^{zP_0}  + z(K_1+D)\otimes
e^{zP_0}  e^{zP_1}P_2 - J\otimes e^{-zP_1}  
\nonumber
\eea
\bea
&&[J,K_i]=\epsilon_{ij}K_j\qquad [J,P_1]=P_2
\qquad [J,P_2]=\frac 1{2z}(e^{-2zP_1}-1)-\frac z2 P_2^2\cr
&&[K_1,P_0]=\frac 1z e^{zP_0}\sinh zP_1 -\frac z2 e^{zP_0}e^{zP_1}P_2^2
\qquad [K_2,P_0]=e^{zP_0}e^{zP_1}P_2\cr
&&[K_1,P_1]=\frac 1z(e^{zP_0}\cosh zP_1-1)-\frac z2 e^{zP_0}e^{zP_1}P_2^2
\qquad [K_1,P_2]=(1-e^{zP_0}e^{-zP_1})P_2\cr
 &&[K_2,P_2]=\frac 1ze^{-zP_1}(e^{zP_0}-\cosh zP_1)+\frac z2 P_2^2
\qquad  [K_2,P_1]=(e^{zP_0}e^{zP_1}-1)P_2\cr
&&[K_1,K_2]=-J\qquad  [D,P_0]=\frac 1z(e^{zP_0}\cosh zP_1-1) + \frac z2
e^{zP_0}e^{zP_1}P_2^2\cr 
&&[D,P_1]=\frac 1z e^{zP_0}\sinh zP_1 + \frac z2
e^{zP_0}e^{zP_1}P_2^2
\qquad [D,P_2]=e^{zP_0}e^{-zP_1}P_2 \cr
&&[P_\mu,P_\nu]=0\qquad [J,P_0]=0 \qquad [D,J]=0 \qquad [D,K_i]=0 .
\label{cb}
\eea

\noindent
b) Special conformal transformations.
\bea
 &&\back\back\Delta(C_0)=1\otimes C_0 + C_0\otimes e^{zP_0}\cosh zP_1 
-C_1\otimes e^{zP_0}\sinh zP_1 - z C_2\otimes e^{zP_0} P_2 \cr
&&  +z(J+K_2)\otimes e^{zP_0}J+  {z^2} (K_1+D)\otimes 
e^{zP_0} e^{zP_1}P_2 J 
 -\frac {z^2}2(C_1-C_0)\otimes 
e^{zP_0} e^{zP_1}P_2^2\cr
&&\back\back\Delta(C_1)=1\otimes C_1 + C_1\otimes e^{zP_0}\cosh zP_1  
-C_0\otimes e^{zP_0}\sinh zP_1   - z C_2\otimes e^{zP_0} P_2 \cr
&&  +z(J+K_2)\otimes e^{zP_0}J+  {z^2} (K_1+D)\otimes 
e^{zP_0} e^{zP_1}P_2 J   -\frac {z^2}2(C_1-C_0)\otimes 
e^{zP_0} e^{zP_1}P_2^2\cr
&&\back\back\Delta(C_2)=1\otimes C_2  + C_2\otimes e^{zP_0}
  + z(C_1-C_0)\otimes
e^{zP_0}  e^{zP_1}P_2  - z(K_1+D)\otimes e^{zP_0}  e^{zP_1} J \cr
&&\qquad  \label{cc}
\eea
\bea
&& [J,C_0]=-zK_1J+\frac z2 J\qquad [J,C_1]=C_2+zDJ\qquad [J,C_2]=-C_1\cr
&& [K_1,C_0]=C_1-\frac z2(K_1+D)+z K_1 D - z (J+K_2)^2\cr
&& [K_2,C_0]=C_2+\frac z2 K_2 + z K_1 J + z(K_1+D) (J+K_2) \cr
&& [K_1,C_1]=C_0-\frac z2  (K_1+D)-\frac z2 (K_1^2 +D^2) 
-\frac z2 (J+K_2)^2\cr
&& [K_2,C_2]=C_0-\frac z2  (K_1+D)-\frac z2 (K_1 +D)^2
 -\frac z2 (J+K_2)^2\cr
&&[K_1,C_2]=z(J+K_2)D\qquad [K_2,C_1]=-zDJ\qquad [P_2,C_2]=-D\cr
&&[P_0,C_0]=D-ze^{zP_0} e^{zP_1}P_2J\qquad 
[P_1,C_1]=-D-ze^{zP_0} e^{zP_1}P_2J\cr
&&[C_1,P_0]=K_1+ z e^{zP_0} e^{zP_1}P_2J\qquad
[C_2,P_0]=K_2-(e^{zP_0} e^{zP_1}-1)J\cr
&&[P_1,C_0]=K_1- z e^{zP_0} e^{zP_1}P_2J\qquad
[P_2,C_1]=-e^{zP_0} e^{-zP_1}J+zDP_2\cr
&&[P_2,C_0]=K_2-zK_1P_2+\frac z2 P_2-(e^{zP_0} e^{-zP_1}-1)J\qquad
[P_1,C_2]=e^{zP_0} e^{zP_1}J\cr
&& [D,C_0]=-C_0+\frac z2    (K_1^2 +D^2) +\frac z2 (J+K_2)^2\label{cd}\\
&&[D,C_1]=-C_1-zK_1D\qquad [D,C_2]=-C_2-z (J+K_2)D\cr
&&[C_1,C_2]=zC_2 - z (J+K_2)C_1+z(K_1+D)C_2\cr
&&[C_0,C_1]=\frac z2 (C_1-C_0) + z (J+K_2)C_2\qquad [C_0,C_2]=
-z(J+K_2)C_1+\frac z2 C_2 .
\nonumber
\eea

Recall that   the  Drinfel'd--Jimbo
$q$-deformation of $so(3,2)$ introduced in  \cite{Lukierskib}   was
performed in a kinematical basis  (as the algebra of the motion group of
the 3+1 anti-de Sitter spacetime) and the two primitive generators were a 
rotation and  the time translation. Now  the primitive generators  are
again  two (the rank of the algebra): the time translation $P_0$ and a
space translation
$P_1$. On the other hand, the symmetry  between $P_\mu$ and
$C_\mu$ is broken in the quantum case (compare to  (\ref{aa})); for
instance, all $P_\mu$ commute among themselves but the $C_\mu$ do not.


\sect{Quantum contractions}

The contraction $U_z({\cal M}_3)\to U_z({\cal G}_3)$  is carried out by
applying the transformations (\ref{ab}) and (\ref{bbd}) to the results
presented in the previous section. Once  the limit $\varepsilon\to 0$ is
taken, the  resultant expressions are rather simplified. The coproduct and
deformed commutation relations of the non-standard quantum conformal
Galilean algebra
$U_z({\cal G}_3)$ are 
\be
\begin{array}{l}
\Delta(X)=1\otimes X+X\otimes 1\qquad \mbox{for}\ X\in\{J,P_\mu,K_i,C_i\}\cr
\Delta(D)=1\otimes D+D\otimes 1  +zK_1\otimes P_1 +z K_2\otimes P_2  \cr
\Delta(C_0)=1\otimes C_0+C_0\otimes 1- zC_1\otimes P_1 -z C_2\otimes P_2  
\end{array}
\label{da}
\ee
\be
[D,P_0]=P_0+\frac z2 (P_1^2+P_2^2)\qquad 
[D,C_0]=-C_0+\frac z2 (K_1^2+K_2^2),
\label{db}
\ee
the remaining commutators are non-deformed and  given by (\ref{ad}). On the
other hand, the element
\be
\begin{array}{l}
{\cal R}=\exp\{r\}=\exp\{z(K_1\wedge P_1 +K_2\wedge P_2)\}\cr
\quad =\exp\{-z P_2\otimes K_2\}\exp\{-z P_1\otimes K_1\}
 \exp\{z K_1\otimes P_1\}\exp\{z K_2\otimes  P_2\}  
\end{array}
\label{dc}
\ee
is a trivial solution of the quantum Yang--Baxter equation since the four
generators involved  commute. Furthermore, it is easy to
check that the property
\be
{\cal R}\Delta (X) {\cal R}^{-1}=\sigma\circ \Delta (X) 
\label{dd}
\ee
is   satisfied for any $X \in U_z({\cal G}_3)$. Then   ${\cal R}$ is a
quantum universal $R$-matrix for 
$U_z({\cal G}_3)$.

Similarly, the contraction $U_z({\cal M}_3)\to U_z({\cal C}_3)$
is provided by the mappings (\ref{ac}) and (\ref{bd}) applied on   the 
conformal Hopf algebra  $U_z({\cal M}_3)$  together the limit $\varepsilon\to
0$. The coproduct and non-vanishing commutation relations of the  quantum
conformal Carroll algebra $U_z({\cal C}_3)$ are given as follows:
\be
\begin{array}{l}
\Delta(X)=1\otimes X+X\otimes 1\qquad \mbox{for}\ X\in\{J,P_\mu\}\cr
\Delta(Y)=1\otimes Y+Y\otimes e^{zP_0}\qquad \mbox{for}\ Y\in\{K_i,C_0\}\cr
\Delta(D)=1\otimes D+D\otimes e^{zP_0} +z K_1\otimes e^{zP_0} P_1 
+ z K_2\otimes e^{zP_0} P_2 \cr
\Delta(C_1)=1\otimes C_1+C_1\otimes e^{zP_0} -z C_0\otimes e^{zP_0} P_1
+z K_2\otimes e^{zP_0} J \cr
\Delta(C_2)=1\otimes C_2+C_2\otimes e^{zP_0} -z C_0\otimes e^{zP_0} P_2 
-z K_1\otimes e^{zP_0} J
\end{array}
\label{de}
\ee
\bea
&&[J,K_i]=\epsilon_{ij}K_j\qquad [J,P_i]=\epsilon_{ij}P_j\qquad 
[J,C_i]=\epsilon_{ij}C_j\cr
&&[K_i,P_i]=\frac 1z(e^{zP_0}-1)\qquad 
[K_i,C_i]=C_0-\frac z2(K_1^2+K_2^2) \qquad [P_0,C_i]=-K_i \cr
&&[C_0,P_i]=-K_i\qquad
[D,P_0]=\frac 1z(e^{zP_0}-1) \qquad [D,P_i]=e^{zP_0}P_i\cr
&&[D,C_0]=-C_0+\frac z2(K_1^2+K_2^2)
\qquad [D,C_i]=-C_i-zK_iD\cr
&&[P_i,C_j]=-\delta_{ij}D+\epsilon_{ij}e^{zP_0}J
\qquad [C_1,C_2]=z(K_1C_2-K_2C_1) .\label{df}
\eea

It is rather remarkable that $U_z({\cal C}_3)$ can be shown  to be
isomorphic to the null-plane quantum Poincar\'e algebra \cite{Null} in the
basis used in
\cite{rnull} by means of
\be
\begin{array}{llll}
P'_+=P_0&\quad P'_-=-C_0&\quad P'_i=K_i&\quad J'_3=-J\cr
K'_3=D&\quad E'_i=-P_i&\quad F'_i=C_i&\quad z'=z/2 
\end{array}
\label{dg}
\ee
where the primed generators and deformation parameter correspond to the 
null-plane quantum Poincar\'e algebra. As a
straightforward consequence the universal
$R$-matrix for $U_z({\cal C}_3)$ (satisfying the quantum  Yang--Baxter
equation and relation (\ref{dd})) reads
\bea
&&{\cal R}=\exp\{-z P_2\otimes K_2\}\exp\{-z P_1\otimes K_1\}
\exp\{- z P_0\otimes D\}\cr
&&\qquad \times \exp\{z D\otimes P_0\}
\exp\{z K_1\otimes P_1\}\exp\{z K_2\otimes  P_2\}. 
\label{dh}
\eea


\sect {Concluding remarks}

Summarizing, we have obtained a new quantum deformation of $so(3,2)$ and we
have related three non-standard quantum conformal algebras via contraction
processes, all of them containing the corresponding Weyl  subalgebra as a
Hopf subalgebra:
\be
U_z( \overline{{\cal G}}(2+1)) \subset U_z({\cal G}_3)
\longleftarrow U_z(\overline{{\cal P}}(2+1) )\subset U_z({\cal M}_3)
\longrightarrow U_z(\overline{{\cal C}}(2+1) )\subset U_z({\cal C}_3) 
\ee
For the contracted quantum algebras the universal $R$-matrices have been
given in a factorized form. The expression of the $R$-matrix associated to
$U_z({\cal M}_3)$ remains as an open problem.

We would like to notice that we could have written the classical
$r$-matrix for $so(2,2)$ (\ref{ba}) as
\be
  r=z(D\wedge P_1 + K_1\wedge P_0);
\ee
indeed, this was exactly the expression chosen   in \cite{Beyond,av}. Its
generalization to $so(3,2)$ would be
\be
r=z(D\wedge P_2+K_1\wedge P_1 + K_2\wedge P_0 + J\wedge P_1).
\label{za}
\ee
From a mathematical point of view, the corresponding quantum deformation is
equivalent to the one  just studied via a simple redefinition of the
generators.  However both deformations exhibit   different physical
features which are stronger when contractions are carried out. More
explicitly, the quantum conformal Galilean and Carroll algebras coming from
(\ref{za}) are no longer equivalent to those above obtained.  For both of
them the transformation of
$z$ would be $z\to \varepsilon^{-2}z$ $(n_0=2)$ leading to the following 
  classical $r$-matrices:
\be 
 {\cal G}_3:\ r=z K_1\wedge P_1  \qquad\quad
  {\cal C}_3:\  r=z  K_2\wedge P_0  .
\ee 
Therefore the latter could not be related to the
null-plane quantum Poincar\'e algebra. The analysis of these and further
possibilities  will be presented elsewhere.



\bigskip \medskip
  
\noindent {\large{{\bf Acknowledgements}}}
 
\bigskip 
\noindent
The author is grateful to Angel Ballesteros for helpful discussions. This
work has been partially supported by DGICYT (Project PB94--1115) from the
Ministerio de Educaci\'on y Ciencia de Espa\~na.

\bigskip 


\end{document}